\documentclass[a4paper, BCOR=0.0mm, DIV=calc]{scrartcl}
\usepackage[T1]{fontenc}
\usepackage[utf8]{inputenc}
\usepackage{url}
\usepackage{lmodern}
\usepackage{microtype}
\usepackage[english]{babel}
\usepackage{amsmath, amssymb, amsfonts}
\usepackage{nicefrac}
\newcommand{\diff}[1]{\mathrm{#1}}
\newcommand{\gam}[1]{\Gamma ( #1 )}

\KOMAoptions{twoside=false, twocolumn=false, headinclude=false, footinclude=false, mpinclude=false, pagesize=auto}
\recalctypearea

\begin{document}
\title{Notes and Remarks on certain logarithmic integrals}
\author{Alexander Aycock}
\date{}
\maketitle
\begin{abstract}
We consider integrals of the form $\int_0^1 \ln{\ln{(\frac{1}{x})}}R{(x)}\diff{d}x$ again, where $R{(x)}$ is a rational function, and we will explain a way to obtain their evaluation. 
\end{abstract}

§1 In mathematics treatises can roughly be divided up into two classes, the first containing those, that expand the boundaries of mathematics and the second containing those, that represent already known things either with new theories in a simpler manner or achieve the same aim with familiar methods in a simpler way.\\

§2 This little memoir falls in the second class, because it considers integrals of the form $\int_0^1 \ln{\ln{(\frac{1}{x}})}R{(x)}\diff{d}x$ again, which have been already been studied by several other mathematicians. Here Vardi \cite{7}, Medina and Moll \cite{17} and Adamchik \cite{18} are to be mentioned, but also a lot of the "oder" mathematicians wrote about such integrals. They mainly treated integrals of the form $\int_0^1\frac{Q{(x)}\diff{d}x}{\ln{x}}$ (where $Q{(x)}$ is a rational function in most cases), which are easily shown to be equivalent to those (at least in certain cases) we want to consider here. We want to mention Euler \cite{11}, Legendre \cite{1} and Malmstèn \cite{3}. The latter contributed the most to this kind of integrals - as we will see below - and actually provided everything to get as far as the mathematicians mentioned earlier. 

So we also want to study this integrals and evaluate them in a simpler manner and want to explain carefully, how to get there a priori.\\

§3 Therefore we will have to look for methods - and will have to explain them -, to obtain these integrals, without using any knowledge from complex theory of functions. It will therefore be convenient to say some things in advance, that will be useful later. These things are, of course, well-known, but the way, to obtain them, that we will present here, seems to be mostly forgotten. And so it will be worth the effort, to describe it.\\

§4 We want to begin with the partial fraction decomposition of the circle or hyperbola functions.

We will start with $\frac{\sin{(ax)}}{\sin{(bx)}}$, the factors of the denominator are $k\pi - bx$ and $k\pi + bx$, where $k$ is a natural number, and it is easily seen, that all factors are simple. Now we put
\[
\frac{\sin{(ax)}}{\sin{(bx)}} = \frac{A}{k\pi - bx} + G(x)
\]
where $A$ is a number, we have to determine, and $G{(x)}$ is a function, not involving the factor $k\pi - bx$. We multiply by $\sin{(bx)}$ and find
\[
\sin{(ax)}=\frac{A\sin{(bx)}}{k\pi-bx}+\sin{(bx)}G{(x)}
\]
Letting $x$ tend to $\frac{k\pi}{b}$ the second term on the right-hand side vanishes; the first can be calculated with L' Hospital's rule, and we find
\[
\sin{(\frac{ak\pi}{b})} = -A\cos{(k\pi)}
\]
Because $k$ is a natural number, we have
\[
A = (-1)^{k+1}\sin{(\frac{ak\pi}{b})}
\]
And in the same way we will find for the other factor $k\pi + bx$, if we set
\[
\frac{\sin({ax})}{\sin({bx})} = \frac{A}{k\pi + bx} + Q(x)
\]
that
\[
B = (-1)^{k+1}\sin{(\frac{ak\pi}{b})}
\]
And by summing over $k$ from $1$ to infinity
\[
\frac{\sin{(ax)}}{\sin{(bx)}} = \sum_{k=1}^{\infty}(-1)^{k+1}\sin{(\frac{ak\pi}{b})}(\frac{1}{k\pi-bx}+\frac{1}{k\pi+bx})
\]
and by contracting the two denominators
\[
\frac{\sin{(ax)}}{\sin{(bx)}} = 2\pi\sum_{k=1}^{\infty}(-1)^{k+1}\sin{(\frac{ak\pi}{b})}(\frac{k}{(k\pi)^2-(bx)^2})
\]
and if we put $x = iy$ we will obtain
\[
\frac{\sinh{(ay)}}{\sinh{(by)}} = \frac{e^{ay}-e^{-ay}}{e^{by}-e^{-by}} = 2\pi\sum_{k=1}^{\infty}(-1)^{k+1}\sin{(\frac{ak\pi}{b})}(\frac{k}{(k\pi)^2+(bx)^2})
\]

§5 We want to consider another example of this kind, $\frac{\cos{(ax)}}{\sin{(bx)}}$, that slightly differs from the preceding one.
Because we see, that the first term of the partial fraction decomposition, because of the zero at $x = 0$ of the denominator and $\cos{(0)} = 1$, will be $\frac{1}{bx}$; for the rest we have exactly the same factors as in the first example. This time we put
\[
\frac{\cos{(ax)}}{\sin{(bx)}} - \frac{1}{bx} = \frac{A}{k\pi - bx} + R(x)
\]
because the first term of the expansion is already known, then we have
\[
\cos{(ax)}-\frac{\sin{(bx)}}{bx}= \frac{A\sin{(bx)}}{k\pi-bx}+\sin{(bx)}R{(x)}
\]
Letting $x$ tend to $\frac{k\pi}{b}$ again we obtain
\[
A = (-1)^{k+1}\cos{(\frac{ak\pi}{b})}
\]
And putting
\[
\frac{\cos{(ax)}}{\sin{(bx)}} - \frac{1}{bx} = \frac{B}{k\pi + bx} + Q(x)
\]
we find
\[
B = -(-1)^{k+1}\cos{(\frac{ak\pi}{b})}
\]
And summing over $k$ from $1$ to infinity again, yields
\[
\frac{\cos{(ax)}}{\sin{(bx)}} = \frac{1}{bx} + \sum_{k=1}^{\infty}(-1)^{k+1}\cos{(\frac{ak\pi}{b})}(\frac{1}{k\pi - bx} - \frac{1}{k\pi + bx})
\]
or after a little simplification
\[
\frac{\cos{(ax)}}{\sin{(bx)}} = \frac{1}{bx} + \sum_{k=1}^{\infty}(-1)^{k+1}\cos{(\frac{ak\pi}{b})}(\frac{2bx}{(k\pi)^2 - (bx)^2})
\]
and the special case $a = b = 1$ gives
\[
\cot{(x)} = \frac{1}{x} + \sum_{k=1}^{\infty}\frac{2x}{x^2-(k\pi)^2}
\]
and in a similar way one can find several more formulas, also for other functions than the circle and hyperbola functions, which subject we will not investigate here. Instead we remark, that, by integrating the last identity and simplifying, we arrive at the famous sine product formula
\[
\sin{(x)} = x\prod_{k=1}^{\infty}(1-\frac{x^2}{(k\pi)^2})
\]
which was given by Euler \cite{9}. For $x = \frac{\pi}{2}$ we find the Wallis product formula for $\frac{\pi}{2}$
\[
\frac{\pi}{2} = \prod_{k=1}^{\infty}(\frac{4k^2}{4k^2-1})
\]
§6 Euler \cite{14} was the first to use the presented method for finding the partial fraction decomposition. At first only for rational functions, later also for the circle and the hyperbola functions. Legendre \cite{1} also arrived at these result with this method. Despite the simplicity of this method it seems to be mostly forgotten, see Sandifer's text \cite{5} for a note on this.\\

§7 We want to add one remark.
Considering the formula
\[
\frac{\sin{(ax)}}{\sin{(bx)}} = 2\pi\sum_{k=1}^{\infty}(-1)^{k+1}\sin{(\frac{ak\pi}{b})}(\frac{k}{(k\pi)^2-(bx)^2})
\]
it can also be interpretated as a function of the variable $a$. So we want to differentiate with respect to $a$, then we will have
\[
x\frac{\cos{(ax)}}{\sin{(bx)}} = \frac{2}{b}\sum_{k=1}^{\infty}(-1)^{k+1}\cos{(\frac{ak\pi}{b})}(\frac{(k\pi)^2}{(k\pi)^2-(bx)^2})
\] 
We already know the value on the left-hand side and writing $(k\pi)^2 = (k\pi)^2-(bx)^2+ (bx)^2$ on the right-hand side, we arrive at the following equation, taking into account the already derived partial fraction decomposition for $\frac{\cos{(ax)}}{\sin{(bx)}}$
\[
x[\frac{1}{bx}+\sum_{k=1}^{\infty}(-1)^{k+1}\cos{(\frac{ak\pi}{b})}(\frac{2bx}{(k\pi)^2 - (bx)^2})]= \frac{2}{b}\sum_{k=1}^{\infty}(-1)^{k+1}\cos{(\frac{ak\pi}{b})}+
\]
\[
\frac{2}{b}\sum_{k=1}^{\infty}(-1)^{k+1}\cos{(\frac{ak\pi}{b})}(\frac{(bx)^2}{(k\pi)^2-(bx)^2})
\]
and this gives after a little simplification
\[
\frac{1}{2} = \sum_{k=1}^{\infty}(-1)^{k+1}\cos{(\frac{ak\pi}{b})}
\]
and for $a = 0$
\[
\frac{1}{2} = \sum_{k=1}^{\infty}(-1)^{k+1}
\]
which series is diverging. But it is well-known, that such series occur very often and in this case we derived it from an identity, which is also known from elsewhere. Therefore we want to see the value of this particular series as $\frac{1}{2}$.\\

§8 It will be convenient to note, that Legendre \cite{1} arrived at the partial fraction decomposition for $\frac{\cos{(ax)}}{\sin{(bx)}}$ using the series $\frac{1}{2} = \sum_{k=1}^{\infty}(-1)^{k+1}$. Euler also explained on several occasions \cite{8}, \cite{16}, how such results should be interpretated and that it makes sense to use such a particular value. These series can be used, if they are interpretated correctly. See Hardy \cite{6} or Ford \cite{19}, to name some more recent mathematicians, who worked on this subject.
We will avoid these series as far as possible, because they require a theory for their explanaition, that cannot be regarded as elementary anymore. But we will nevertheless use the value of the one series, that we found here, later. \\

§9 Before we are able to get to our main results, we have to show some identities, that will be useful later. \\

We start with the formulae 
\[
\frac{\sin{(a)}}{1+2y\cos{(a)}+y^2} = \sum_{n=1}^{\infty}(-1)^{n-1}y^{n-1}\sin{(na)}
\]
and
\[
\frac{a}{2}= \sum_{n=1}^{\infty}(-1)^{n-1}\frac{\sin{(na)}}{n}
\]
which are well-known in the theory of Fourier series, but they can be proven by elementary means and the proof, we will give here, traces back to Euler \cite{16}. Consider the geometric series
\[
\frac{1}{1+x} =   \sum_{n=0}^{\infty}(-1)^{n}x^{n}
\]
And for $x = ye^{i\phi}$ by using the Euler identity $e^{ix}=\cos{(x)}+i\sin{(x)}$ this gives
\[
\frac{1}{1+y\cos{(\phi)}+iy\sin{(\phi)}} = \sum_{n=0}^{\infty}(-1)^{n}y^{n}e^{in\phi}
\]
By using de Moivre's formula and after a simplification
\[
\frac{1+y\cos{(\phi)}}{1+2y\cos{(\phi)}+y^2}-\frac{i\sin{(\phi)}}{1+2y\cos{(\phi)}+y^2} = \sum_{n=0}^{\infty}(-1)^{n}y^{n}(\cos{(n\phi)}+i\sin{(n\phi)})
\]
Comparing the real and imaginary parts yields
\[
\frac{1+y\cos{(\phi)}}{1+2y\cos{(\phi)}+y^2}= \sum_{n=0}^{\infty}(-1)^{n}y^{n}\cos{(n\phi)}
\]
and
\[
\frac{\sin{(\phi)}}{1+2y\cos{(\phi)}+y^2}=\sum_{n=1}^{\infty}(-1)^{n-1}y^{n-1}\sin{(n\phi)}
\]
The first formula gives the divergent series from above, if we put $y = 1$, and the second is the one, we wanted to demonstrate.
To get the formula for $\frac{a}{2}$ we can either integrate the forst formula with respect to $\phi$ or the second with respect to $y$ and put $y = 1$ after the integration.
We omit the exact calculation here, because it the one given should be sufficient.
We could also derive the identity in exactly the same way as the first from the power series for $\ln{(1+x)}$. But this would lead to the question, how we find the logarithm of a complex number, already requiring a little complex analysis. Therefore we will not do it that way, but mention, that Euler \cite{16} actually did it this way.\\

§10 Another useful identity for evaluating logarithmic integrals is this one
\[
\int_0^{\infty}\frac{e^{au}-e^{-au}}{e^{\pi u}-e^{-\pi u}}\cos{(uz)}\diff{d}u= \frac{\sin{(a)}e^{-z}}{1+2e^{-z}\cos{(a)}+e^{-2z}}
\]
To show this one, we will need the following integral for positive $n$
\[
\int_0^{\infty}\frac{\cos{(nx)}\diff{d}x}{1+x^2}= \frac{\pi}{2}e^{-n}
\]
We only have to cnsider it for positve $n$, because the function is symmetric in $n$ and we will just need it for positive $n$ later anyway.
We suppose, that we do not know the answer yet, but we know, the value will be a function of $n$, and so we set
\[
\int_0^{\infty}\frac{\cos{(nx)}\diff{d}x}{1+x^2}= f(n)
\]
then we will also have by differentiating
\[
\int_0^{\infty}\frac{x\sin{(nx)}\diff{d}x}{1+x^2}= -f^{\prime}(n)
\]
By integrating by parts one easily confirms the following fomulas
\[
\int_0^{\infty}e^{-at}\cos{(bt)}\diff{d}t = \frac{a}{a^2+b^2}
\]
and 
\[
\int_0^{\infty}e^{-at}\sin{(bt)}\diff{d}t = \frac{b}{a^2+b^2}
\]
And we find as a special case
\[
\int_0^{\infty}e^{-xt}\sin{(t)}\diff{d}t = \frac{1}{1+x^2}
\]
Multiplying this equation by $\cos{(nx)}$ and integrating from $0$ to infinity with respect to $x$ yields 
\[
\int_0^{\infty}\int_0^{\infty}\cos{(nx)}e^{-xt}\sin{(t)}\diff{d}t\diff{d}x=\int_0^{\infty}\frac{\cos{(nx)}\diff{d}x}{1+x^2}
\]
The right-hand side is $f{(n)}$ amd therefore we have
\[
f(n) = \int_0^{\infty}(\int_0^{\infty}\cos{(nx)}e^{-xt}\diff{d}x)\sin{(t)}\diff{d}t
\]
The inner integral can be expressed with the formula above. Using this, this leads to
\[
f(n)=\int_0^{\infty}\frac{t\sin{(t)}\diff{d}t}{n^2+t^2}
\]
and for $t =ny$
\[
f(n)= \int_0^{\infty}\frac{y\sin{(ny)}\diff{d}y}{1+y^2}=f^{\prime}(n)
\]
So we arrived at a differential equation for $f{(n)}$, whose solution is seen to be $Ce^{-n}$, where $C$ is a constant, but we have
\[
Ce^{-0}= C = \int_0^{\infty}\frac{\diff{d}y}{1+y^2} = \frac{\pi}{2}
\]
which we wanted to prove, because now we have
\[
\int_0^{\infty}\frac{\cos{(nx)}\diff{d}x}{1+x^2}= \frac{\pi}{2}e^{-n}
\]
§11 Even if this proof was not completely rigorous and some conditions are to be added, it leads to the right result, which was already known to mathematician like Legendre \cite{1}, who also evaluated it without the use of complex function theory.
But this result is - of course - easyly confirmed by the calculus of residues nowadays. \\

§12 But now we are able, to prove our important identity, from which essentially all other things will flow; this will also seen below. We want to demonstrate now
\[
\int_0^{\infty}\frac{e^{ax}-e^{-ax}}{e^{\pi x}-e^{-\pi x}}\cos{(nx)}\diff{d}x= \frac{\sin{(a)}e^{-n}}{1+2e^{-n}\cos{(a)}+e^{-2n}}
\]
We start from the integral and attempt to evaluate it, we found the partial fraction for $\frac{e^{ax}-e^{-ax}}{e^{\pi x}-e^{-\pi x}}$ above, if we put $b= \pi$ in that formula. Replacing this expression with its partial fraction decomposition gives
\[
\int_0^{\infty}\frac{e^{ax}-e^{-ax}}{e^{\pi x}-e^{-\pi x}}\cos{(nx)}\diff{d}x= \frac{2}{\pi}\sum_{k=1}^{\infty}(-1)^{k+1}\sin{(ak)}k\int_0^{\infty}\frac{\cos{(nx)}\diff{d}x}{k^2+x^2}
\]
Setting $yk =x$ in the intgral on the right-hand side 
\[
\int_0^{\infty}\frac{e^{ax}-e^{-ax}}{e^{\pi x}-e^{-\pi x}}\cos{(nx)}\diff{d}x= \frac{2}{\pi}\sum_{k=1}^{\infty}(-1)^{k+1}\sin{(ak)}\int_0^{\infty}\frac{\cos{(nky)}\diff{d}y}{1+y^2}
\]
This integral can be evaluated and we will have
\[
\int_0^{\infty}\frac{e^{ax}-e^{-ax}}{e^{\pi x}-e^{-\pi x}}\cos{(nx)}\diff{d}x= \frac{2}{\pi}\sum_{k=1}^{\infty}(-1)^{k+1}\sin{(ak)}\frac{\pi}{2}e^{-nk}
\]
The arising series can be expressed in finite terms using the results from above
and yields
\[
\int_0^{\infty}\frac{e^{ax}-e^{-ax}}{e^{\pi x}-e^{-\pi x}}\cos{(nx)}\diff{d}x= \frac{\sin{(a)}e^{-n}}{1+2e^{-n}\cos{(a)}+e^{-2n}}
\]
This is the one, we wanted to show and it will be observed later, that it will lead to the disired evaluation of logarithmic integrals, how Malmstèn \cite{3} realized at first, even if this partiular formula traces back to Legendre \cite{1}. \\

§13 In the papers of Vardi \cite{7}, Medina and Moll \cite{17} and Adamchik \cite{18} one will see quite fast, that the logarithmic integrals follow from functiona equations of certain Dirichlet series, series of the form $\sum_{k=1}^{\infty}\frac{a{(n)}}{n^{s}}$. So, if we show these, we can claim, to have everything in our, what the mentioned authors proved, concerning the evaluation at least. We just would have to follow their way. Now, so I claim, we can already derive these functional equations with few sketches from the preceding. It will therefore be worth the effort to show these functiona equations, before explaining another related method.\\

§14 For this purpose we note, that we have 
\[
\int_0^{\infty}z^{s-1}\cos{(uz)}\diff{d}z= \frac{\Gamma{(s)}}{u^{s}}\cos{(\frac{s\pi}{2})}
\]
This can be seen as follows. We have, where - as usual - $\Gamma{(s)}$ is the well-known function, defined by the integral $\int_0^{\infty}e^{-t}t^{s-1}\diff{d}t $, if we write $kt$ for $t$ in this integral,
\[
\frac{\Gamma{(s)}}{k^{s}}=\int_0^{\infty}e^{-kt}t^{s-1}\diff{d}t
\]
and therefore for $k = a+bi$ and using the Euler identity again
\[
\frac{\Gamma{(s)}}{(a+bi)^{s}}=\int_0^{\infty}e^{-at}\cos{(bt)}t^{s-1}\diff{d}t -i\int_0^{\infty}e^{-at}\sin{(bt)}t^{s-1}\diff{d}t
\]
And we have, of course, 
\[
\frac{\Gamma{(s)}}{(a+bi)^{s}}\cdot\frac{(a-bi)^{s}}{(a-bi)^{s}}=\frac{\Gamma{(s)}}{(a^2+b^2)^{s}}(a-bi)^{s}
\]
writing the complex number in polar coordinates this reduces to
\[
\frac{\Gamma{(s)}}{(a+bi)^{s}}= \frac{\Gamma{(s)}}{(a^2+b^2)^{\frac{s}{2}}}(\cos{(s\arctan{\frac{b}{a}})}-i\sin{(s\arctan{\frac{b}{a}})})
\]
which is the value of our integral; if we compare the real parts, we will obtain the desired identity, after having changed the letters
\[
\int_0^{\infty}e^{-xz}z^{s-1}\cos{(uz)}\diff{d}z= \frac{\Gamma{(s)}}{(x^2+u^2)^{\frac{s}{2}}}\cos{(s\arctan{\frac{u}{x}})}
\]
Letting $x$ tend to $0$, the argument of the arctan tends to infinity and therefore the arctan tends to $\frac{\pi}{2}$. This shows our formula
\[
\int_0^{\infty}z^{s-1}\cos{(uz)}\diff{d}z= \frac{\Gamma{(s)}}{u^{s}}\cos{(\frac{s\pi}{2})}
\]
which goes back to Euler \cite{12} again and gives the Fresnel integrals for $s=\frac{1}{2}$.
We now want to multiply the last equation by $\frac{e^{au}-e^{-au}}{e^{\pi u}-e^{-\pi u}}$ and integrate from $0$ to infintity with respect to $u$, which yields
\[
\int_0^{\infty}\int_0^{\infty}\frac{e^{au}-e^{-au}}{e^{\pi u}-e^{-\pi u}}z^{s-1}\cos{(uz)}\diff{d}z\diff{d}u= \Gamma{(s)}\cos{(\frac{s\pi}{2})}\int_0^{\infty}\frac{e^{au}-e^{-au}}{e^{\pi u}-e^{-\pi u}}\frac{\diff{d}u}{u^{s}}
\]
The double integral can be simplified by using the idendity from above, namely,
\[
\int_0^{\infty}\frac{e^{au}-e^{-au}}{e^{\pi u}-e^{-\pi u}}\cos{(uz)}\diff{d}u= \frac{\sin{(a)}e^{-z}}{1+2e^{-z}\cos{(a)}+e^{-2z}}
\]
and gives
\[
\int_0^{\infty}\frac{e^{au}-e^{-au}}{e^{\pi u}-e^{-\pi u}}\frac{\diff{d}u}{u^{s}}=\frac{1}{\Gamma{(s)}\cos{(\frac{s\pi}{2})}}\int_0^{\infty}\frac{z^{s-1}\sin{(a)}e^{-z}\diff{d}z}{1+2e^{-z}\cos{(a)}+e^{-2z}}
\]
For $e^{-z}= \ln{(y)}$ this reduces to
\[
\int_0^{\infty}\frac{e^{au}-e^{-au}}{e^{\pi u}-e^{-\pi u}}\frac{\diff{d}u}{u^{s}}=
\frac{1}{\Gamma{(s)}\cos{(\frac{s\pi}{2})}}\int_0^{1}\frac{\ln^{s-1}{(\frac{1}{y})}\sin{(a)}\diff{d}y}{1+2y\cos{(a)}+y^{2}}
\]
\\

§15 And this is already the fundmaental formula for the functional equations of certain Dirichlet series, which are important for the consideration about logarithmic integrals; this can be seen as follows. We have
\[
\int_0^{\infty}\frac{e^{au}-e^{-au}}{e^{\pi u}-e^{-\pi u}}\frac{\diff{d}u}{u^{s}}=
\frac{\sin{(a)}}{\Gamma{(s)}\cos{(\frac{s\pi}{2})}}\int_0^{1}\frac{\ln^{s-1}{(\frac{1}{y})}\diff{d}y}{1+2y\cos{(a)}+y^{2}}
\]
and
\[
\int_0^{\infty}\frac{e^{-\pi u}(e^{au}-e^{-au})}{1-e^{-2\pi u}}\frac{\diff{d}u}{u^{s}}=
\frac{\sin{(a)}}{\Gamma{(s)}\cos{(\frac{s\pi}{2})}}\int_0^{1}\frac{\ln^{s-1}{(\frac{1}{y})}\diff{d}y}{1+2y\cos{(a)}+y^{2}}
\]
If we use the series expansions for $\frac{1}{1-e^{-2\pi u}}$, and $\frac{\sin{(a)}}{1+2y\cos{(a)}+y^{2}}$ - we found the series for this second expression above -, we find
\[
\sum_{n=0}^{\infty}\int_0^{\infty}(e^{-2n\pi u-\pi u -a}-e^{-2n\pi u -\pi u  +a})u^{-s}\diff{d}u= \frac{1}{\Gamma{(s)}\cos{(\frac{s\pi}{2})}}\sum_{n=1}^{\infty}(-1)^{n-1}\sin{(na)}\int_0^{1}\ln^{s-1}{(\frac{1}{y})} y^{n-1} \diff{d}y
\]
All occuring integrals can be expressed in terms of $\Gamma{(s)}$ and this leads to
\[
\Gamma{(1-s)}\sum_{n=0}^{\infty}(\frac{1}{((2n+1)\pi-a)^{1-s}}-\frac{1}{((2n+1)\pi+a)^{1-s}})=\frac{\Gamma{(s)}}{\Gamma{(s)}\cos{(\frac{s\pi}{2})}}\sum_{n=1}^{\infty}(-1)^{n-1}\frac{\sin{(na)}}{n^{s}}
\]
On the right-hand side the $\Gamma{(s)}$ functions cancel, and writing $\pi -a$ instead of $a$, this gives a similar formula 
\[
\Gamma{(1-s)}\sum_{n=0}^{\infty}(\frac{1}{(2n\pi+a)^{1-s}}-\frac{1}{((2n+2)\pi-a)^{1-s}})=\frac{1}{\cos{(\frac{s\pi}{2})}}\sum_{n=1}^{\infty}\frac{\sin{(na)}}{n^{s}}
\]
And finally putting $a =\frac{k\pi}{l}$ and writing $1-s$ instead $s$ we arrive a these two formulas
\[
\sum_{n=0}^{\infty}(\frac{1}{((2n+1)l-k)^{s}}-\frac{1}{((2n+1)l+k)^{s}})= \frac{(\frac{\pi}{l})^{s}}{\sin{(\frac{s\pi}{2})}\Gamma{(s)}}\sum_{n=1}^{\infty}(-1)^{n+1}\frac{\sin{(\frac{nk\pi}{l})}}{n^{1-s}}
\]
and
\[
\sum_{n=0}^{\infty}(\frac{1}{(2nl+k)^{s}}-\frac{1}{((2n+2)l+k)^{s}})= \frac{(\frac{\pi}{l})^{s}}{\sin{(\frac{s\pi}{2})}\Gamma{(s)}}\sum_{n=1}^{\infty}\frac{\sin{(\frac{nk\pi}{l})}}{n^{1-s}}
\]
§16 If we use the first series for an example and put $k=1$, $l=2$ and $k=1$, $l=3$ respectively, we obtain these to equations
\[
\sum_{n=0}^{\infty}\frac{(-1)^{n}}{(2n+1)^{s}}= \frac{(\frac{\pi}{2})^{s}}{\sin{(\frac{s\pi}{2})}\Gamma{(s)}}\sum_{n=0}^{\infty}\frac{(-1)^{n}}{(2n+1)^{1-s}}
\]
and
\[
\sum_{n=0}^{\infty}(\frac{1}{(3n+1)^{s}}-\frac{1}{(3n+2)^{s}})= \frac{(\frac{2\pi}{3})^{s}}{\sin{(\frac{s\pi}{2})}\Gamma{(s)}}\sum_{n=0}^{\infty}(\frac{1}{(3n+1)^{1-s}}-\frac{1}{(3n+2)^{1-s}})
\]
These and some more were at first given by Malmstèn \cite{2} and proved rigorously. So Malmstèn was probably the first person to proove functional equation of Dirichlet series, which cannot be praised enough.\\\\

§17 It will be observed, that you can add many more, even infinitely many, to those two, we gave here, and certainly all, that were considered by Vardi \cite{7}, Medina and Moll \cite{17} and Adamchik \cite{18}. Only the functional equation of the famous Riemann zeta function, $\zeta{(s)}=\sum_{n=1}^{\infty}\frac{1}{n^{s}}$, or equivalently the Dirichlet eta function, $\eta{(s)}=\sum_{n=1}^{\infty}\frac{(-1)^{n+1}}{n^{s}}$, are missing, both being such examples, that cannot be obtained from our formulas by putting in values.\\

§18 But you can perform a little trick, as I demonstrated on another occasion, considering Malmstèn's paper \cite{2}. I will nevertheless, because it fits right in, be convenient, to show this little trick again and show
\[
\eta{(1-s)}=\frac{2^{s}-1}{1-2^{s-1}}\pi^{-s}\cos{(\frac{s\pi}{2})}\Gamma{(s)}\eta{(s)}
\]
For the sake of brevity we want to set
\[
\lambda{(s)}=\sum_{n=0}^{\infty}\frac{1}{(2n+1)^{s}}=\frac{2^{s}-1}{2^{s}-2}\eta{(s)}
\]
We will start from our fundamental identity again
\[
\int_0^{\infty}\frac{e^{au}-e^{-au}}{e^{\pi u}-e^{-\pi u}}\frac{\diff{d}u}{u^{s}}=
\frac{1}{\Gamma{(s)}\cos{(\frac{s\pi}{2})}}\int_0^{1}\frac{\ln^{s-1}{(\frac{1}{y})}\sin{(a)}\diff{d}y}{1+2y\cos{(a)}+y^{2}}
\]
We want to divide both sides by $\sin{(a)}$ and let $a$ tend to $0$, where we use the limit
\[
\lim_{a \rightarrow 0}\frac{e^{au}-e^{-au}}{\sin{a}} = 2u
\]
which can be derived from L' Hospitals's rule. This yields 
\[
2\cos{\left(\frac{\pi s}{2}\right)}\int_0^{\infty}\frac{u^{1-s}\diff{d}u}{e^{\pi u}-e^{-\pi u}} = \frac{1}{\gam{s}}\int_0^1\frac{\ln^{s-1}{(\frac{1}{y})}\diff{d}y}{(1+y)^2}
\]
Multiplying the numerator and the denominator by $e^{-\pi u}$ and using the series expansions for $\frac{1}{1-e^{-2\pi u}}$ and $\frac{1}{(1+y)^{2}}$ we find
\[
2\cos{\left(\frac{\pi s}{2}\right)}\int_0^{\infty}\sum_{n=0}^{\infty}e^{-(2n+1)\pi u}u^{1-s}\diff{d}u = \frac{1}{\gam{s}}\int_0^1\ln^{s-1}{\left(\frac{1}{y}\right)\sum_{n=1}^{\infty}(-1)^{n+1}n y^{n-1}\diff{d}y}
\]
Interchanging the order of summation and integration, both integrals can be reduced to constant multiples of $\Gamma{(s)}$ and yield
\[
\frac{2\cos{\left(\frac{\pi s}{2}\right)}}{\pi^{2-s}}\sum_{n=0}^{\infty}\frac{1}{(2n+1)^{2-s}}{\gam{2-s}} 
= \frac{1}{\gam{s}}\sum_{n=1}^{\infty}(-1)^{n+1}n\frac{\gam{s}}{n^s}\;.
\]
And with our abbreviations
\[
\frac{2\cos{\left(\frac{\pi s}{2}\right)}}{\pi^{2-s}}\lambda{(2-s)}\gam{2-s} = \eta{(s-1)}
\]
If we finally write $2-s$ instead of $s$ and replace $\lambda{(s)}$ with the identity involving $\eta{(s)}$, we will get
\[
\eta{(1-s)} = \frac{2^s-1}{1-2^{s-1}}\pi^{-s}\cos{\left(\frac{\pi s}{2}\right)}\gam{s}\eta{(s)}
\]
which equation we wanted to show. We want to note after all this, that this functional equation and the one for $\beta{(s)}= \sum_{n=0}^{\infty}\frac{(-1)^{n}}{(2n+1)^{s}}$ were already known by Euler \cite{15} and he found it with the help of divergent series. Further you can show by replacing $\eta{(s)}$ with the relation to $\zeta{(s)}$ - we have $\zeta{(s)}= \frac{\eta{(s)}}{1-2^{1-s}}$ , of course - the functional equation for $\zeta{(s)}$, namely
\[
\zeta{(1-s)}=\frac{2}{(2\pi)^{s}}\Gamma{(s)}\cos{(\frac{s\pi}{2})}\zeta{(s)}
\]
which was at first proved rigorously by Riemann in his famous memoir \cite{2}.\\

§19 Now we see, that we have everything in our hands, that Vardi \cite{7}, Medina and Moll \cite{17} and Adamchik \cite{18} needed, to obtain their results, we only representes everything with integrals, so that we could end the paper right here, having provided a a priori method to evaluate the logarithmic integrals.\\

§20 But these integrals can also be calculated in a different way, as special values of the Fourier series expansion for $\ln{\Gamma{(x)}}$.
The Fourier series was given by Kummer \cite{4} in 1847, but also by Malmstèn \cite{3} one year earlier in an entirely different way, which we will essentially present here, on the one hand, because it involves logarithmic integrals - and therefore reveals the connection of these intgrals to the Fourier series - and on the other hand, because it deserves to be mentionend, containing some interesting manipulations. We go on to the proof.\\

§21 We want to derive the Fourier series for $\ln{\frac{\Gamma{(\frac{1}{2}+x)}}{\Gamma{(\frac{1}{2}-x)}}}$. For this we note the following identities, which are straight-forward to prove, see also Legendre \cite{1} and Malmstèn \cite{3}.
\[
\int_0^{1}\frac{y^{a}-y^{b}}{1-y^2}\diff{d}y = \frac{1}{2}[\psi{(\frac{1}{2}(b+1))}-\psi{(\frac{1}{2}(a+1))}]
\]
where $\psi{(x)}$ is the Digamma function and can be defined by
\[
\psi{(x)}= -\gamma +\int_0^{1}\frac{1-t^{x-1}}{1-t}\diff{d}t= \frac{\diff{d}}{\diff{d}x}\ln{\Gamma{(x)}}
\]
From this we find
\[
\psi{(\frac{1}{2})}=-\gamma -2\ln{2}= -\gamma -\ln{4}
\]
and we also have
\[
\int_0^{1}\ln{\ln{(\frac{1}{y})}}y^{n-1}\diff{d}y= -\frac{\gamma+\ln{n}}{n}
\]
which can be shown by differentiating the following formula with respect to $s$ and setting $s=1$ afterwards
\[
\int_0^{1}\ln^{s-1}{(\frac{1}{y})}y^{n-1}\diff{d}y= \frac{\Gamma{(s)}}{n^s}
\]
where $\gamma$ is, as above, $-\Gamma^{\prime}{(1)}$. Actually, of course, $\gamma$ is the famous Euler-Mascheroni constant and is defined by the following limit
\[
\gamma = \lim_{n \rightarrow \infty}(\sum_{k=1}^{n}\frac{1}{k}-\ln{n}) = 0,577215664901...
\]
and it can be shown, that we also have $\gamma= -\Gamma^{\prime}{(1)}$. But, because we will neither use the exact value nor the limit definition, and only the value $-\Gamma^{\prime}{(1)}$, we can omit this particular proof.\\

§22 Having said all this in advance, we can derive the Fourier series expansion.
From the preceeding we already have
\[
\cos{(\frac{s\pi}{2})}\int_0^{\infty}\frac{e^{au}-e^{-au}-2au}{e^{\pi u}-e^{-\pi u}}\frac{\diff{d}u}{u^{s}}=\frac{1}{\Gamma{(s)}}\int_0^{1}(\frac{\ln^{s-1}{(\frac{1}{y})}\sin{(a)}}{1+2y\cos{(a)}+y^{2}}-\frac{a\ln^{s-1}{(\frac{1}{y})}}{(1+y)^2})\diff{d}y
\]
We differentiate with respect to $s$ and obtain
\[
-\frac{\pi}{2}\sin{(\frac{s\pi}{2})}\int_0^{\infty}\frac{e^{au}-e^{-au}-2au}{e^{\pi u}-e^{-\pi u}}\frac{\diff{d}u}{u^{s}}-\cos{(\frac{s\pi}{2})}\int_0^{\infty}\frac{e^{au}-e^{-au}-2au}{e^{\pi u}-e^{-\pi u}}\frac{\ln{u}\diff{d}u}{u^{s}}
\]
\[
=-\frac{\Gamma^{\prime}{(s)}}{\Gamma{(s)}}\int_0^{1}(\frac{\ln^{s-1}{(\frac{1}{y})}\sin{(a)}}{1+2y\cos{(a)}+y^{2}}-\frac{a\ln^{s-1}{(\frac{1}{y})}}{(1+y)^2})\diff{d}y+
\]
\[
\frac{1}{\Gamma{(s)}}\int_0^{1}\ln{\ln{(\frac{1}{y})}}\ln^{s-1}{(\frac{1}{y})}(\frac{\sin{(a)}}{1+2y\cos{(a)}+y^{2}}-\frac{a}{(1+y)^2})\diff{d}y
\]
Because the second integral on the richt-hand side of the equation keeps bounded as $s$ tends to $1$, we get, letting $s$ tend to $1$
\[
-\frac{\pi}{2}\int_0^{\infty}\frac{e^{au}-e^{-au}-2au}{e^{\pi u}-e^{-\pi u}}\frac{\diff{d}u}{u}=\gamma\int_0^{1}(\frac{\sin{(a)}}{1+2y\cos{(a)}+y^{2}}-\frac{a}{(1+y)^2})\diff{d}y
\]
\[
+\int_0^{1}\frac{\ln{\ln{(\frac{1}{y})}}\sin{(a)}\diff{d}y}{1+2y\cos{(a)}+y^{2}}-a\int_0^{1}\frac{\ln{\ln{(\frac{1}{y})}}\diff{d}y}{(1+y)^{2}}
\]
The first integral on the right-hand side of the equation is seen to vanish by direct calculation, for the second we want for the sake of brevity write $F{(a)}$.
Then our equation becomes
\[
-\frac{\pi}{2}\int_0^{a}\frac{\diff{d}}{\diff{d}a}(\int_0^{\infty}\frac{e^{au}-e^{-au}-2au}{e^{\pi u}-e^{-\pi u}}\frac{\diff{d}u}{u})\diff{d}a=-a\int_0^{1}\frac{\ln{\ln{(\frac{1}{y})}}\diff{d}y}{(1+y)^{2}}+F(a)
\]
And if we expand $\frac{1}{(1+y)^2}$ into a series
\[
-\frac{\pi}{2}\int_0^{a}\int_0^{\infty}\frac{e^{au}-e^{-au}-2}{e^{\pi u}-e^{-\pi u}}\diff{d}u\diff{d}a= a\sum_{n=1}^{\infty}(-1)^{n}n\int_{0}^{1}\ln{\ln{(\frac{1}{y})}}y^{n-1}\diff{d}y+ F(a)
\]
If we set $e^{-\pi u} =y$ and write $y^0+y^0$ instead of $2$ and evaluate the integral on the right-hand side, with the identity mentioned earlier, we will have
\[
-\frac{1}{2}\int_0^{a}\int_0^{1}\frac{y^{\frac{a}{\pi}}+y^{-\frac{a}{\pi}}-y^0-y^0}{1-y^2}\diff{d}y\diff{d}a=a\sum_{n=1}^{\infty}(-1)^{n+1}n\frac{\gamma+\ln{n}}{n}
\]
The inner integral, depending on $y$, can be expressed in terms of $\psi{(x)}$ using the identity from above, the occuring sum $\sum_{n=1}^{\infty}(-1)^{n+1}$ is equal to $\frac{1}{2}$, as we saw earlier.
Therefore
\[
-\frac{1}{2}\int_0^{a}\frac{1}{2}(2\psi{(\frac{1}{2})}-\psi{(\frac{1}{2}+\frac{a}{2\pi})}-\psi{(\frac{1}{2}-\frac{a}{2\pi})})\diff{d}a=\frac{\gamma a}{2}+\frac{a}{2}\sum_{n=1}^{\infty}(-1)^{n+1}\ln{n}+F(a)
\]
The occuring sum, involving logarithms, leads to this infinite product
\[
\prod_{k=1}^{\infty}(\frac{2k-1}{2k})
\]
But it can be reduces to the Wallis product formula as follows
\[
\prod_{k=1}^{\infty}(\frac{2k-1}{2k})=\sqrt{\prod_{k=1}^{\infty}(\frac{4k^2}{4k^2-1})}=\sqrt{\frac{2}{\pi}}
\]
Simplyfying and using these values gives
\[
-\frac{1}{4}\cdot 2\psi{(\frac{1}{2})}a+\frac{2\pi}{4}\ln{\Gamma{(\frac{1}{2}+\frac{a}{2\pi})}}-\frac{2\pi}{4}\ln{\Gamma{(\frac{1}{2}-\frac{a}{2\pi})}}=\frac{\gamma a}{2}+\frac{a}{2}\ln{\sqrt{\frac{2}{\pi}}}+F(a)
\]
Solving for $F{(a)}$ and subsituting the value for $\psi{(\frac{1}{2})}$, yields after some simplification
\[
F(a)=\int_0^{1}\frac{\ln{\ln{(\frac{1}{y})}}\sin{(a)}\diff{d}y}{1+2y\cos{(a)}+y^{2}}=\frac{\pi}{2}\ln{(\frac{(2\pi)^{\frac{a}{\pi}}\ln{\Gamma{(\frac{1}{2}+\frac{a}{2\pi})}}}{\ln{\Gamma{(\frac{1}{2}-\frac{a}{2\pi})}}})}
\]
If we apply the series for $\frac{\sin{(a)}\diff{d}y}{1+2y\cos{(a)}+y^2}$, we obtain
\[
\sum_{n=1}^{\infty}(-1)^{n-1}\sin{(na)}\int_0^{1}\ln{\ln{(\frac{1}{y})}}y^{n-1}\diff{d}y=\frac{\pi}{2}\ln{(\frac{(2\pi)^{\frac{a}{\pi}}\ln{\Gamma{(\frac{1}{2}+\frac{a}{2\pi})}}}{\ln{\Gamma{(\frac{1}{2}-\frac{a}{2\pi})}}})}=F(a)
\]
The integral is known again and leads to
\[
F(a)=-\gamma\sum_{n=1}^{\infty}(-1)^{n-1}\frac{\sin{(na)}}{n}-\sum_{n=1}^{\infty}(-1)^{n-1}\frac{\sin{(na)}\ln{n}}{n}
\]
We saw the first series to be $\frac{a}{2}$. If we use this and let $a =2\pi x$ and solve for $\ln{\frac{\Gamma{(\frac{1}{2}+x)}}{\Gamma{(\frac{1}{2}-x)}}}$, we arrive at the following Fourier series expansion
\[
\ln{\frac{\Gamma{(\frac{1}{2}+x)}}{\Gamma{(\frac{1}{2}-x)}}}=-2x(\gamma+\ln{(2\pi)})+\frac{2}{\pi}\sum_{k=1}^{\infty}(-1)^{k}\frac{\ln{k}}{k}\sin{(2k\pi x)}
\]
And this series enables us, to evaluate a lot of logarithmic integrals, because we know see the connection between the integrals and the series from all the things we derived. So it will be convenient to offer some examples, that it can be seen better.\\

§23 The first example should be Vardi's integral, as the following integral is now called,
\[
\int_0^{1}\frac{\ln{\ln{(\frac{1}{y})}}\diff{d}y}{1+y^{2}}=\frac{\pi}{2}\ln{(\frac{(2\pi)^{\frac{1}{2}}\ln{\Gamma{(\frac{3}{4})}}}{\ln{\Gamma{(\frac{1}{4})}}})}
\]
We can use the integral for the evaluation
\[
F(a)=\int_0^{1}\frac{\ln{\ln{(\frac{1}{y})}}\sin{(a)}\diff{d}y}{1+2y\cos{(a)}+y^{2}}=\frac{\pi}{2}\ln{(\frac{(2\pi)^{\frac{a}{\pi}}\ln{\Gamma{(\frac{1}{2}+\frac{a}{2\pi})}}}{\ln{\Gamma{(\frac{1}{2}-\frac{a}{2\pi})}}})}
\]
we get formula 1 in Vardi's paper, because $a=\frac{\pi}{2}$ already leads directly to the desired integral.
To get more formulas, we just have to put in more values for $a$, we want to take $a=\frac{\pi}{3}$, this will give
\[
\int_0^{1}\frac{\ln{\ln{(\frac{1}{y})}}\diff{d}y}{1+y+y^{2}}=\frac{\pi}{\sqrt{3}}\ln{(\frac{(2\pi)^{\frac{1}{3}}\ln{\Gamma{(\frac{2}{3})}}}{\ln{\Gamma{(\frac{1}{3})}}})}
\]
and for $a=\frac{2\pi}{3}$ we will get this formula
\[
\int_0^{1}\frac{\ln{\ln{(\frac{1}{y})}}\diff{d}y}{1-y+y^{2}}=\frac{\pi}{\sqrt{3}}\ln{(\frac{(2\pi)^{\frac{2}{3}}\ln{\Gamma{(\frac{7}{6})}}}{\ln{\Gamma{(-\frac{1}{6})}}})}
\]
And this, with the known properties of $\ln{\Gamma{(x)}}$, reduces to the following value, also given by Vardi in his paper \cite{7}.
\[
\int_0^{1}\frac{\ln{\ln{(\frac{1}{y})}}\diff{d}y}{1+y+y^{2}}=\frac{2\pi}{\sqrt{3}}[\frac{5}{6}\ln{(2\pi)}-\ln{\Gamma(\frac{1}{6})}]
\] 
We would for example use the reflection formula for $\Gamma{(x)}$, namely
\[
\Gamma{(x)}\Gamma{(1-x)}=\frac{\pi}{\sin{\pi x}}
\]
which was given by Euler \cite{13}, which we did not show and therefore will not use here. \\

§24 It will be convenient, to add one more example, that involves a rational function, whose denominator is raised to the second power.
For this purpose we consider our integral again
\[
F(a)=\int_0^{1}\frac{\ln{\ln{(\frac{1}{y})}}\sin{(a)}\diff{d}y}{1+2y\cos{(a)}+y^{2}}=\frac{\pi}{2}\ln{(\frac{(2\pi)^{\frac{a}{\pi}}\ln{\Gamma{(\frac{1}{2}+\frac{a}{2\pi})}}}{\ln{\Gamma{(\frac{1}{2}-\frac{a}{2\pi})}}})}
\]
If we differentiate it with respect to $a$ again, we obtain
\[
\int_0^{1}\frac{\ln{\ln{(\frac{1}{y})}}((1+y^2)\cos{(a)}+2y)\diff{d}y}{(1+2y\cos{(a)}+y^{2})^2}= \frac{\ln{(2\pi)}}{2}+\frac{1}{4}\psi{(\frac{1}{2}+\frac{a}{2\pi})}+\frac{1}{4}\psi{(\frac{1}{2}-\frac{a}{2\pi})}
\]
and if we put $a=\frac{\pi}{2}$ and use the known values (or just calculate them from the integral expression for $psi{(x)}$)
\[
\psi{(\frac{1}{4})}=-\gamma-\frac{\pi}{2}-\ln{8}
\] 
and
\[
\psi{(\frac{3}{4})}=-\gamma+\frac{\pi}{2}-\ln{8}
\] 
we get this formula
\[
\int_0^{1}\frac{y\ln{\ln{(\frac{1}{y})}}\diff{d}y}{(1+y)^2}=\frac{1}{2}\ln{\sqrt{\pi}}-\frac{\ln{4}}{4}-\frac{\gamma}{4}
\]
which formula is a special case of a more general one, given by Adamchik in propostion 5 in his paper \cite{18}.\\

§25 The first three integrals were all given by Vardi \cite{7}, the first two also by Malmstèn \cite{3}, including Vardi's integral. And from his formulas there follow many others.\\

§26 It is now easyly seen, that we provided everything again, to get many logarithmic integrals. The greatest problem consists in the evaluation of certain sums, involving logarithms, which - as we saw - can be expressed in finite terms with the given Fourier series expansion.\\

We could also evaluate integrals as this one
\[
\int_0^{1}\frac{y\ln{\ln{(\frac{1}{y})}}\diff{d}y}{(1-y+y^2)^2}=-\frac{\gamma}{3}-\ln{(\frac{6\sqrt{3}}{\pi})}+\frac{\pi \sqrt{3}}{27}[5\ln{(2\pi)}-6\ln{\Gamma(\frac{1}{6})}]
\]
if we allowed the use of divergent series or at least the principle of analytic continuation, because we would be lead to certain divergent series, which correspond to certain sums and quotients of geometric series and their derivatives at the point $x=1$.\\

§27 But having shown at least a little bit about these integrals and having proved the functional equations for some Dirichlet series, we will put aside more concrete evaluations, because they are well-presented by Medina and Moll \cite{17} and Adamchik \cite{18} and we would need divergent series, which require a deeper theory, that cannot be regarded as elementary anymore and that is not as straight-forward as the calculations, we did and explained, in this memoir. 
We will give more evaluations on another occasion and use Euler's definition \cite{8} of a divergent series and its sum and will see, whether it leads to the same results.\\

{\bf Acknowledgement:} I would like to thank Sebastian Koch and Arseny Skryagin for showing interest in this little paper and making many useful suggestions. I would also like to thank Artur Diener for helping a lot, writing this paper.

\end{document}